\documentclass[twoside,leqno,10pt]{article}

\usepackage{graphics,ifthen}
\usepackage{latexsym}
\usepackage{amsmath}
\usepackage{amssymb}
\usepackage{mathrsfs}
\usepackage{natbib}
\begin{document}
\input{latexP.sty}
\input{referencesP.sty}
\input {epsf.sty}
\def\ind{\stackrel{\mathrm{ind}}{\sim}}
\def\iid{\stackrel{\mathrm{iid}}{\sim}}
\def\Prodi{\mathop{{\lower9pt\hbox{\epsfxsize=15pt\epsfbox{pi.ps}}}}}
\def\prodi{\mathop{{\lower3pt\hbox{\epsfxsize=7pt\epsfbox{pi.ps}}}}}
\def\Definition{\stepcounter{definitionN}
    \Demo{Definition\hskip\smallindent\thedefinitionN}}
\def\EndDefinition{\EndDemo}
\def\Example#1{\Demo{Example [{\rm #1}]}}
\def\EndExample{\qed\EndDemo}
\def\Category#1{\centerline{\Heading #1}\rm}
%% Paper specific definitions
\def\e{\text{\hskip1.5pt e}}
\newcommand{\eps}{\epsilon}
\newcommand{\proof}{\noindent {\bf Proof:\ }}
\newcommand{\remarks}{\noindent {\bf Remarks:\ }}
\newcommand{\note}{\noindent {\bf Note:\ }}
\newcommand{\examp}{\noindent {\bf Example:\ }}
\newcommand{\Lower}[2]{\smash{\lower #1 \hbox{#2}}}
\newcommand{\ben}{\begin{enumerate}}
\newcommand{\een}{\end{enumerate}}
\newcommand{\bi}{\begin{itemize}}
\newcommand{\ei}{\end{itemize}}
\newcommand{\hp}{\hspace{.2in}}
\newcommand{\DoubleRLarrow}[1]{\Lower{-0.01in}{$\underleftarrow{\Lower{0.07in}{$\overrightarrow{\vspace*{0.15in}\hspace*{#1}}$}}$}}
\newcommand{\DoubleLRarrow}[1]{\Lower{-0.01in}{$\underrightarrow{\Lower{0.07in}{$\overleftarrow{\vspace*{0.15in}\hspace*{#1}}$}}$}}
\newcommand{\SingleRarrow}[1]{\Lower{-0.08in}{$\underrightarrow{\hspace*{#1}}$}}
\newcommand{\SingleLarrow}[1]{\Lower{-0.08in}{$\underleftarrow{\hspace*{#1}}$}}

%Mathrsfs Font
\newcommand{\Bcr}{\mathscr{B}}
\newcommand{\Ucr}{\mathscr{U}}
\newcommand{\Gcr}{\mathscr{G}}
\newcommand{\Dcr}{\mathscr{D}}
\newcommand{\CS}{\mathscr{C}}
\newcommand{\Fcr}{\mathscr{F}}
\newcommand{\Icr}{\mathscr{I}}
\newcommand{\Lcr}{\mathscr{L}}
\newcommand{\Mcr}{\mathscr{M}}
\newcommand{\Ncr}{\mathscr{N}}
\newcommand{\Pcr}{\mathscr{P}}
\newcommand{\Qcr}{\mathscr{Q}}
\newcommand{\Scr}{\mathscr{S}}
\newcommand{\Tcr}{\mathscr{T}}
\newcommand{\Xcr}{\mathscr{X}}
\newcommand{\Wcr}{\mathscr{W}}
%Mathbb Font
\newcommand{\E}{\mathbb{E}}
\newcommand{\F}{\mathbb{F}}
\newcommand{\I}{\mathbb{I}}
\newcommand{\Q}{\mathbb{Q}}
\newcommand{\X}{\mathbb{X}}
\newcommand{\Pe}{\mathbb{P}}
\newcommand{\M}{\mathbb{M}}
\newcommand{\R}{\mathbb{R}}
\newcommand{\Wbb}{\mathbb{W}}

%Math abbreviations
\newcommand{\Xn}{X_1,\ldots,X_n}
\newcommand{\lan}{\langle}
\newcommand{\ra}{\rangle}
\newcommand{\edr}{\mathrm{e}}
\newcommand{\ddr}{\mathrm{d}}
\newcommand{\Var}{\operatornamewithlimits{Var}}
\newcommand{\ld}{\ldots}
\newcommand{\cd}{\cdots}
\def\simind{\stackrel{\mbox{\scriptsize{ind}}}{\sim}}
\def\simiid{\stackrel{\mbox{\scriptsize{iid}}}{\sim}}

%Greek letters
\newcommand{\Ga}{\Gamma}
\newcommand{\si}{\sigma}
\newcommand{\la}{\lambda}
\newcommand{\ph}{\varphi}
\newcommand{\ep}{\varepsilon}
\newcommand{\ga}{\gamma}
\newcommand{\Om}{\Omega}
\newcommand{\om}{\omega}

\newtheorem{thm}{Theorem}[section]
\newtheorem{defin}{Definition}[section]
\newtheorem{prop}{Proposition}[section]
\newtheorem{lem}{Lemma}[section]
\newtheorem{cor}{Corollary}[section]
\newcommand{\rb}[1]{\raisebox{1.5ex}[0pt]{#1}}
\newcommand{\mc}{\multicolumn}
\def\Beta{\text{Beta}}
\def\Dir{\text{Dirichlet}}
\def\DP{\text{DP}}
\def\P{{\bf p}}
\def\fhat{\widehat{f}}
\def\GA{\text{gamma}}
\def\ind{\stackrel{\mathrm{ind}}{\sim}}
\def\iid{\stackrel{\mathrm{iid}}{\sim}}
\def\J{{\bf J}}
\def\K{{\bf K}}
\def\min{\text{min}}
\def\N{\text{N}}
\def\p{{\bf p}}
\def\Sj{{\bf S}}
\def\sj{{\bf s}}
\def\U{{\bf U}}
\def\u{{\bf u}}
\def\w{{\bf w}}
\def\W{{\bf W}}
\def\X{{\bf X}}
\def\x{{\bf x}}
\def\y{{\bf y}}
\def\Y{{\bf Y}}
\newcommand{\reals}{{\rm I\!R}}
\newcommand{\PR}{{\rm I\!P}}
\def\Z{{\bf Z}}
\def\yy{{\mathcal Y}}
\def\rr{{\mathcal R}}
\def\mm{{\mathcal M}}
\def\BP{\text{beta}}
\def\ts{\tilde{t}}
\def\Ns{\tilde{N}}
\def\Ps{\tilde{P}}
\def\gs{\tilde{g}}
\def\fs{\tilde{f}}
\def\ys{\tilde{Y}}
\def\ps{\tilde{P}}
\def
\Report{\centerline{\small{\rm PD COAG-FRAG}}}
\def\Author{\centerline{\small{\rm M-W. HO, L.F. JAMES, J.W. LAU}}} \pagestyle{myheadings} \markboth{\Author}{\Report}
\thispagestyle{empty} \bct\Heading Coagulation Fragmentation Laws
Induced By General Coagulations of Two-Parameter Poisson-Dirichlet
Processes\lbk\lbk\smc {\sc Man-Wai Ho, Lancelot F. James, and John
W. Lau}\footnote{ 
\indent\eightit Keywords and phrases. \rm
         Cauchy-Stieltjes transforms,
         coagulation-fragmentation,
         exchangeable random partitions,
         Poisson-Kingman models,
         random discrete distributions,
         two-parameter Poisson Dirichlet
         }
\lbk\lbk \BigSlant National University of Singapore, Hong Kong
University of Science and Technology
and University of Bristol\rm \lbk %(\today)
\ect \Quote Pitman~(1999) describes a duality relationship between
fragmentation and coagulation operators. An explicit relationship is
described for the two-parameter Poisson-Dirichlet laws, say
{\footnotesize $PD(a,b)$}, with parameters {\footnotesize
$(\alpha,\theta)$} and {\footnotesize $(\beta,\theta/\alpha)$},
wherein {\footnotesize $PD(\alpha, \theta)$} is coagulated by
{\footnotesize $PD(\beta,\theta/\alpha)$} for {\footnotesize
$0<\alpha<1$}, {\footnotesize $0 \leq\beta<1$} and {\footnotesize
$-\beta<\theta/\alpha$.} This remarkable explicit agreement was
obtained by combinatorial methods via exchangeable partition
probability functions~(EPPF). It has been noted that such a method
is not easy to employ for more general processes. This work
discusses an alternative analysis which can feasibly extend the
characterizations above to more general models of {\footnotesize
$PD(\alpha,\theta)$} coagulated with some law {\footnotesize $Q$}.
The analysis exploits distributional relationships between
compositions of species sampling random probability measures and
coagulation operators. It is shown, based on results of Vershik, Yor
and Tsilevich~(2004) and James~(2002), how the calculation of
generalized Cauchy-Stieltjes transforms of random probability
measures provides a blueprint to obtain explicit characterizations
of {\footnotesize $PD(\alpha,\theta)$} coagulated with some law
{\footnotesize $Q$}. We use this to obtain explicit descriptions in
the case where {\footnotesize $Q$} corresponds to a large class of
power tempered Poisson Kingman models described in James~(2002).
That is, explicit results are obtained for models outside of the
{\footnotesize $PD(\beta,\theta/\alpha)$} family. We obtain a new
proof of Pitman's result as a by-product. Furthermore, noting an
obvious distinction from the class of {\footnotesize
$PD(\alpha,\theta)$} derived from a stable subordinator, we discuss
briefly the case of Dirichlet processes coagulated by various
{\footnotesize $Q$}.
 \EndQuote

\section{Introduction}
 Let ${\mathcal P}^{\downarrow}_{1}=\{p=(p_{i}):p_{1}\ge
p_{2}\ge p_{_3}\ldots\ge 0; \sum_{i=1}^{\infty}p_{i}=1\}$.
Furthermore, for a sequence $(x_{1},x_{2},\ldots)$ of non-negative
real numbers with $\sum_{i=1}^{\infty}x_{i}=1$, let
$RANK(x_{1},x_{2},\ldots)\in {\mathcal P}^{\downarrow}_{1}$ be the
decreasing rearrangement of terms of the sequence. Pitman~(1999,
2005), in particular section 5.4 of Pitman~(2005), gives the
following definition of coagulation and fragmentation kernels on
${\mathcal P}^{\downarrow}_{1}$. For each probability measure $Q$ on
${\mathcal P}^{\downarrow}_{1}$, two Markov-transition kernels
$Q-COAG$ and $Q-FRAG$ can be defined for ${\mathcal
P}^{\downarrow}_{1}$ as follows. For $p\in {\mathcal
P}^{\downarrow}_{1}$, $(Q-COAG)(p,\cdot)$ is the distribution on
${\mathcal P}^{\downarrow}_{1}$ of $RANK(\sum_{i}p_{i}I(U_{i}\in
I^{Q}_{j}),j\geq 1)$ where $(I^{Q}_{j})$ is a $Q$-partition of
$[0,1]$ and the $U_{i}$ are i.i.d. uniform on $[0,1]$ independent of
$(I^{Q}_{j}).$ Note that for brevity we refer the reader to
Pitman~(2005, ch.5) for further explanations of the above
quantities. $(Q-FRAG)(p,\cdot)$ is the distribution of
$RANK(p_{i}Q_{ij}, i, j\ge 1)$ where $(Q_{ij})_{j\ge 1}$ has
distribution $Q$ for each $i$, and these sequences are independent
as $i$ varies.

For a probability measure $P$ on ${\mathcal P}^{\downarrow}_{1}$,
let $R:=P(Q-COAG)$. That is the random probability measure on
${\mathcal P}^{\downarrow}_{1}$ defined by
$$
R(\cdot)=\int_{{\mathcal P}^{\downarrow}_{1}}P(dr)Q-COAG(r,\cdot)
$$
which may be called $P$ \emph{coagulated} by $Q$. In principle,
there are many ways to characterize the laws $R$, $P$, $Q-COAG$,
$Q-FRAG$. For example, one may do this via their corresponding
exchangeable partition probability functions~(EPPF) on the space
of partitions of the integers. However this is, in general, a
non-trivial matter. Ideally one wants to identify such laws which
may be described nicely via the graph of Pitman~(2005),
\begin{eqnarray}
  & \tilde{Q} &  \nonumber\\
  X & \DoubleRLarrow{0.5in} & Y \label{arrow}\\
  & {\hat Q} &\nonumber
\end{eqnarray}
where one reads~\mref{arrow} as $\Pe(Y\in \cdot|X)={\tilde
Q}(X,\cdot)$ and $\Pe(X\in \cdot|Y)={\hat Q}(Y,\cdot)$. With
respect to the present context, $X$ has distribution $P$, ${\tilde
Q}(X,\cdot)=Q-COAG(X,\cdot)$, $Y$ has distribution $R$ and ${\hat
Q}(Y,\cdot)=Q-FRAG(Y,\cdot)$. The general task, that we shall
consider, is given $P$ and $Q-COAG$, find $R$ and $Q-FRAG.$

Pitman~(1999) establishes the most general known
coagaluation/fragmentation duality of this type using the
two-parameter Poisson Dirichlet distribution on ${\mathcal
P}^{\downarrow}_{1}.$ The two-parameter Poisson Dirichlet
distribution, denoted as $PD(\alpha,\theta)$, for the separate
ranges $0\leq \alpha<1$, $\theta> -\alpha$ and $\alpha=-\kappa$,
$\theta=m\kappa$ for $\kappa>0$ and some integer $m=1,2,\ldots$ is
discussed in for instance Pitman and Yor~(1997) and Pitman~(2005)
and has numerous applications and interpretations. We shall
provide more details shortly. First Theoerm 12 of Pitman~(1999)
may be described in terms of the following diagram as given in
Pitman~(2005); for $0<\alpha<1, 0\leq \beta<1,
-\beta<\theta/\alpha$,
\begin{eqnarray}
  & PD(\beta,\theta/\alpha)-COAG &  \nonumber\\
  PD(\alpha,\theta) & \DoubleRLarrow{0.5in} & PD(\alpha\beta,\theta) \label{Pitman}\\
  & PD(\alpha,-\alpha\beta)-FRAG &\nonumber
\end{eqnarray}
where the notation $PD(\alpha,\theta)$ and $PD(\alpha\beta,\theta)$
in~\mref{Pitman} is to be understood as some $X$ and $Y$ having
these respective laws. In other words, this gives the explicit
description of the coagulation/fragmentation duality in relation to
$PD(\alpha,\theta)$ coagulated by $PD(\beta,\theta/\alpha).$ That is
one can set, $P=PD(\alpha,\theta)$, $Q=PD(\beta,\theta/\alpha)$ and
$R=PD(\alpha\beta,\theta)$. Pitman~(1999) proves this result via a
combinatorial argument involving the respective EPPF's of the
various two-parameter Poisson-Dirichlet models. The argument used
exploited the Gibbs structure of these EPPF's and as noted by
Pitman~(1999) is not obviously extendable to obtain explicit
expressions for other laws.

In this paper we show how one may replace the combinatorial
argument by an argument involving generalized Cauchy-Stieltjes
transform and moreover extend Pitman's result to more general
models of $PD(\alpha,\theta)$ coagulated by some $Q$. That is to
say given $Q$ we want to complete the description of the following
diagram
\begin{eqnarray}
  & Q-COAG &  \nonumber\\
  PD(\alpha,\theta) & \DoubleRLarrow{0.5in} & Y\label{arrow3}.\\
  & Q-FRAG &\nonumber
\end{eqnarray}

A key to our exposition is the following characterization via
exchangeable random probability measures. First every random
sequence $(P_{i})\in {\mathcal P}^{\downarrow}_{1}$ has a law $P$
which determines and is determined by the law of the random
probability measure
$$
\tau_{P}(\cdot)=\sum_{i=1}^{\infty}P_{i}\delta_{U_{i}}(\cdot)
$$
where $U_{i}$ are iid uniform $[0,1]$. This representation is
equivalent to saying that $\tau_{P}$ is a \emph{species sampling
random probability measure}~[see Pitman~(1996)] based on a Uniform
distribution. Now associating the definition of random probability
measure $\tau_{Q}$ with $Q$ in an obvious way, Lemma 5.18 of
Pitman~(2005) states that the law $R=P(Q-COAG)$ is the unique
probability distribution on ${\mathcal P}^{\downarrow}_{1}$ such
that
$$
(\tau_{R}(u), 0\leq u\leq
1)\overset{d}=(\tau_{P}(\tau_{Q}(u)),0\leq u\leq 1)
$$
where it is assumed that $(\tau_{P}(u),0\leq u\leq 1)$ and
$(\tau_{Q}(u),0\leq u\leq 1)$ are independent. We shall also use the
notation $\tau_{R}=\tau_{P}\circ \tau_{Q}$ to denote composition.
See also Bertoin and Pitman~(2000), Bertoin and Le Gall~(2003, 2005)
for a related discussion.

Our approach is to try to ascertain directly the distribution of
$\tau_{R}=\tau_{P}\circ \tau_{Q}$, when $\tau_{P}$ is determined
by a $PD(\alpha,\theta)$ model. The main tool will be the explicit
evaluation of the generalized Cauchy-Stieltjes transform for
$\tau_{R}$. As we shall show, this approach is particularly well
suited for the $PD(\alpha,\theta)$ models due to results of
Vershik, Yor and Tsilevich~(2004) model in conjunction with the
results of James~(2002). We apply the Cauchy-Stieltjes transforms
in James~(2002) to extend~\mref{Pitman} to a family of $Q$
belonging to a class of power tempered Poisson Kingman laws. This
constitutes a large class of models which are derived from rather
arbitrary continuous infinitely divisible random variables. As an
important example, we show that when $Q$ is a Dirichlet process,
Pitman's result in~\mref{Pitman} for $Q=PD(0,\theta/\alpha)$
follows from the identity of Cifarelli and Regazzini~(1990) and a
new more general result for $Q=PD(0,\nu)$, $\nu>\theta/\alpha,$
follows by a characterization of the Dirichlet process given in
James~(2005).

Some other notable, but not exhaustive, list of references for
various types of coagulation/fragmentation models include Aldous and
Pitman~(1998), Bolthausen and Sznitman~(1998), Bertoin~(2002),
Bertoin and Goldschmidt~(2004), Dong, Goldschmidt and Martin~(2005)
and Schweinsberg~(2000).

\section{Construction of Poisson Kingman Type Random Probability Measures}
We first describe the class of models $Q$ we shall explicitly
consider. Let $T$ denote a strictly positive random variable with
density denoted as $f_{T}$ and Laplace transform
$$
\E[{\mbox e}^{-\lambda T}]={\mbox
e}^{-\psi(\lambda)}=\int_{0}^{\infty}{\mbox
e}^{-t\lambda}f_{T}(t)dt
$$
where $\psi(\lambda)=\int_{0}^{\infty}(1-{\mbox e}^{-\lambda
s})\rho(ds)$ and $\rho$ denotes its unique L\'evy density. Let
$H(\cdot)$ denote a probability measure on a Polish space $\Xcr$.
For the moment we shall assume that $H$ is fixed and non-atomic. We
will later relax this assumption. It is known that for each $T$ and
fixed $H$ one may construct a finite completely random measure, say
$\mu,$ on a Polish space $\Xcr,$ characterized by its Laplace
functional for every positive measureable function $g$ on $\Xcr$ as
\Eq \E[{\mbox e}^{-\mu(g)}|H]={\mbox
e}^{-\int_{\Xcr}\psi(g(x))H(dx)} \label{laplacem}\EndEq where
$\mu(g)=\int_{\Xcr}g(x)\mu(dx).$ It is evident that
$T=\mu(\Xcr):=\int_{\Xcr}\mu(dx)=\int_{\Xcr}I\{x\in \Xcr\}\mu(dx)$.
We denote the law of $\mu$ as $\Pe(d\mu|\rho H)$, where
$\Pe(\cdot|\rho H)$ is a probability measure on a suitably
measureable space of finite measures, say $\Mcr.$ Harkening back to
Kingman~(1975) one may describe a class of random probability
measures on $\Xcr$ by the normalization \Eq
P_{K}(\cdot)=\frac{\mu(\cdot)}{T}=\sum_{i=1}^{\infty}P_{i}\delta_{Z_{i}}(\cdot)
\label{PKgen} \EndEq where $(P_{i})\in {\mathcal
P}^{\downarrow}_{1}$ has some law denoted as $Q=PK(\rho)$ and
independent of $P_{i}$ the $(Z_{i})$ are iid $H$.  That is to say
the $P_{K}$ constitute a class of species sampling random
probability models. The construction of the $(P_{i})$ equates with
the \emph{basic} Poisson-Kingman models discussed in Pitman~(2003).
Pitman~(2003) provides a thorough characterization of the laws
$Q=PK(\rho)$ on ${\mathcal P}^{\downarrow}_{1}$ via their
corresponding exchangeable partition probability function~(EPPF).
Specifically, according to Corollary 6 of Pitman~(2003),  for some
random partition of the integers ${1,\ldots,n}$, $(A_{1},\ldots,
A_{k})$, with block sizes $|A_{i}|=n_{i}$ for $i=1,\ldots, k\leq n$
blocks, the EPPF associated with each $Q$ is given by
$$
p_{K}(n_{1},\ldots, n_{k}):=
\frac{(-1)^{n-k}}{\Gamma(n)}\int_{0}^{\infty}\lambda^{n-1}{\mbox
e}^{-\psi(\lambda)}\prod_{i=1}^{k}\psi_{n_{i}}(\lambda)d\lambda
$$
where for $m=1,\ldots,n,$ \Eq
\psi_{m}(\lambda):=\frac{d^m}{d\lambda^m}\psi(\lambda) =
(-1)^{m-1}\kappa_{m}(\lambda), \label{psim} \EndEq and,
$$
\kappa_{m}(\lambda)=\int_{0}^{\infty}s^{m}{\mbox e}^{-\lambda
s}\rho(ds)
$$
represents the $m$-th cumulant of a random variable with
\emph{tilted} density ${\mbox e}^{\psi(\lambda)} {\mbox
e}^{-\lambda t}f_{T}(t).$ As discussed in Pitman~(1996), the EPPF
$p_{K}$ along with  specific knowledge of $H$ determines the law
of $P_{K}$ which is governed by $\Pe(d\mu|\rho H).$ The basic
Poisson-Kingman laws generate a much larger class of laws by first
conditioning $(P_{i})|T=t$ or equivalently $P_{K}|T=t$ and
substituting $f_{T}(t)dt$ by another probability measure on
$(0,\infty)$, $\gamma(dt)$. Pitman~(2003) denotes these laws as
$PK(\rho,\gamma)=\int_{0}^{\infty}PK(\rho|t)\gamma(dt)$. The
$PK(\rho,\gamma)$ is referred to as a \emph{Poisson-Kingman
distribution with L\'evy density $\rho$ and mixing distribution}
$\gamma.$

\Remark It is obvious that if the $(Z_{i})$ are replaced by
$(U_{i})$ then the composition of such a
$P_{K}(\cdot)=\sum_{i=1}^{\infty}P_{i}\delta_{U_{i}}(\cdot)$
random probability measure with an $H$, $P_{K}\circ H$, is
equivalent to a random probability measure determined by
$\int_{\Xcr}\psi(g((x))H(dx)$ as above. Equivalently, any $P_{K}$
in~\mref{PKgen}, can be represented as
$$
\sum_{k=1}^{\infty}P_{i}\delta_{Z_{i}}(\cdot)\overset
{d}=\frac{\tilde{\mu}{(H(\cdot))}}{\tilde{\mu}((H(\Xcr))}
$$
where the law of ${\tilde \mu}$ on $[0,1]$ is specified by its
Laplace functional with $\E[{\mbox e}^{-\tilde \mu(g)}]:={\mbox
e}^{-\int_{0}^{1}\psi(g((x))dx}.$ That is $\mu={\tilde \mu}\circ
H,$ with distribution characterized by~\mref{laplacem}.
Importantly these result hold for any fixed $H$, whether it
possesses atoms or not. \EndRemark
\subsection{Two-parameter Poisson-Dirichlet models}
The $PD(\alpha,\theta)$ models for $0\leq\alpha<1$ and
$\theta>-\alpha$ are special cases of the above construction, that
is $PK(\rho,\gamma)$ models. The two-parameter $(\alpha,\theta)$
Poisson-Dirichlet random probability measure, with parameters
$0\leq \alpha<1$ and $\theta> -\alpha$ has the known
representation,
$$
P_{\alpha,\theta}(dx)=
\frac{\mu_{\alpha,\theta}(dx)}{T_{\alpha,\theta}}
$$
where $\mu_{\alpha,\theta}$ is a finite random measure on $\Xcr$
with law denoted as $\Pe_{\alpha,\theta}(\cdot|H)$, and
$T_{\alpha,\theta}=\mu_{\alpha,\theta}(\Xcr)$ is a random variable.
The law of the random measure $\mu_{\alpha,\theta}$ can be described
as follows. When $\alpha=0$, $\mu_{0,\theta}$ is a Gamma process
with shape $\theta H$, hence $P_{0,\theta}$ is a Dirichlet process
with shape $\theta H$. That is a $PD(0,\theta)$ model for $(P_{i})$
coupled with a specification for $H$ yields a Dirichlet process with
shape parameter $\theta H$, for $\theta>0$. In this case of
$\mu_{0,\theta},$ the $\psi(g(x))$ is expressed for any positive
measureable function $g$ as \Eq d_{\theta}(g(x)):=\theta
\ln(1+g(x))\label{DPpsi},\EndEq and its L\'evy density is
$\rho_{0,\theta}(ds)=\theta s^{-1}{\mbox e}^{-s}ds$. Hence its law
is $\Pe_{0,\theta}(\cdot|H):=\Pe(\cdot|\rho_{0,\theta}H).$ The total
random mass, say $T_{0,\theta}=\mu_{0,\theta}(\Xcr),$ is a Gamma
random variable with shape parameter $\theta.$ For the
$PD(\alpha,0)$ model, recall that
$$
\rho_{\alpha,0}(ds)=\frac{\alpha
s^{-\alpha-1}}{\Gamma(1-\alpha)}s^{-\alpha-1}ds
$$
is the L\'evy density corresponding to a stable law of index
$0<\alpha<1$. Equivalently, $T_{\alpha,0}$ is a stable random
variable determined by $\rho_{\alpha,0}$, with density denoted as
$f_{\alpha}(t)=f_{T_{\alpha,0}}(t)$, its Laplace transform is given
by
$$
\E[{\mbox e}^{-\lambda T_{\alpha,0}}]={\mbox
e}^{-\lambda^{\alpha}} .$$ This shows that $\mu_{\alpha,0}$ is a
completely random measure based on a stable law.  The law of
$\mu_{\alpha,0}$ is
$\Pe_{\alpha,0}(\cdot|H):=\Pe(\cdot|\rho_{\alpha,0}H).$ In the
cases above both $\mu_{\alpha,0}$ and $\mu_{0,\theta}$ are
completely random measures and $PD(\alpha,0)=PK(\rho_{\alpha,0})$
and $PD(0,\theta):=PK(\rho_{0,\theta}).$ This is not the case for
$PD(\alpha,\theta)$ models for the range $0<\alpha<1$ and
$\theta\neq 0$, $\theta>-\alpha$.  The two-parameter
Poisson-Dirichlet model with $0<\alpha<1$ and $\theta\neq 0$,
$\theta>-\alpha$ is obtained by the specification
$PD(\alpha,\theta):=PK(\rho_{\alpha,0},f_{\alpha,\theta})$ where
$$
f_{\alpha,\theta}(t)dt=c_{\alpha,\theta}t^{-\theta}f_{\alpha}(t)dt,
$$
with $
c_{\alpha,\theta}=1/\E[T^{-\theta}_{\alpha,0}]=\Gamma(\theta+1)/\Gamma(\frac{\theta}{\alpha}+1).$
At the level of the random measure $\mu_{\alpha,\theta}$ one has
the following absolute continuity relationship, for every
measureable function $h$, \Eq
\E[h(\mu_{\alpha,\theta})|H]=c_{\alpha,\theta}\E[T^{-\theta}_{\alpha,0}h(\mu_{\alpha,0})|H]
=c_{\alpha,\theta}\int_{\Mcr}T^{-\theta}h(\mu)\Pe(d\mu|\rho_{\alpha,0}H)
\label{absPD} \EndEq where the first expectation is taken with
respect to the law $\Pe_{\alpha,\theta}(\cdot|H).$ The
$PD(\alpha,\theta)$ model is defined in general for two ranges
$0\leq \alpha<1$ and $\theta>-\alpha$ or $\alpha=-\kappa<0$, and
$\theta=m\kappa$ for $m=1,2,\ldots.$ In any case the EPPF is given
by $$
p_{\alpha,\theta}(n_{1},\ldots,n_{k})=\frac{(\theta+\alpha)_{k-1\uparrow\alpha}
\prod_{i=1}^{k}(1-\alpha)_{n_{i}-1\uparrow1}}{{(\theta+1)}_{n-1\uparrow1}}
$$ where $(x)_{n\uparrow \alpha}:=\prod_{i=0}^{n-1}(x+i\alpha).$

\subsection{$PK(\rho,\gamma_{\theta})$ models} James~(2002, section 6),
influenced by the power tempering construction of the
$PD(\alpha,\theta)$ models, discussed and analyzed various
features of a natural extension to more general $PK(\rho,\gamma)$
models of this type. Suppose that for a $PK(\rho)$ model there
exists $-\infty<\theta<\infty$ such that
$$
\frac{1}{m_{\theta}(\rho)}=\int_{0}^{\infty}t^{-\theta}f_{T}(t)dt<\infty,
$$
then one may define a class of power tempered PK models by
specifying $PK(\rho,\gamma_{\theta})$ with
$$\gamma_{\theta}(dt)= m_{\theta}(\rho)t^{-\theta}f_{T}(t)dt.$$ Hereafter, we
shall only consider the range $\theta>-\alpha.$ The corresponding
random probability measures on $\Xcr$ are denoted as\Eq
P_{K,\theta}(\cdot)=\frac{\mu(\cdot)}{T}=\sum_{i=1}^{\infty}P_{i}\delta_{Z_{i}}
(\cdot) \label{PKT} \EndEq where $(P_{i})\in{\mathcal
P}^{\downarrow}_{1}$ has law $PK(\rho,\gamma_{\theta})$ and
$(Z_{i})$ are iid $H.$ Equivalently if $\Pe(d\mu|\rho H)$ denotes
the distribution of $\mu$ under the the $PK(\rho)$ model, for a
specific $H$, then we say that
$$
\Pe(d\mu|\rho H,\gamma_{\theta})=m_{\theta}(\rho)
T^{-\theta}\Pe(d\mu|\rho H)
$$
is the distribution of $\mu$ under the $PK(\rho,\gamma_{\theta})$
laws. That is, similar to~\mref{absPD}, if $\mu_{K,\theta}$ denotes
a version of $\mu$ with law $\Pe(\cdot|\rho H,\gamma_{\theta})$,
then one has the following absolute continuity relationship, for
every measureable function $h$, \Eq
\E[h(\mu_{K,\theta})|H]=m_{\theta}(\rho)\E[T^{-\theta}h(\mu)|H]
=m_{\theta}(\rho)\int_{\Mcr}T^{-\theta}h(\mu)\Pe(d\mu|\rho H).
\label{absPK} \EndEq Hence we see that~\mref{absPK} is a
generalization of~\mref{absPD}. It follows from Pitman~(2003) that
the EPPF of these models may be described as
$$
p_{K,\theta}(n_{1},\ldots,
n_{k}):=\frac{(-1)^{n-k}m_{\theta}(\rho)}{\Gamma(\theta+n)}\int_{0}^{\infty}\lambda^{\theta+n-1}{\mbox
e}^{-\psi(\lambda)}\prod_{i=1}^{k}\psi_{n_{i}}(\lambda)d\lambda,
$$
where $\psi_{n_{i}}(\lambda)$ is defined as in~\mref{psim}.
\section{Cauchy-Stieltjes Transforms} We now proceed to show
how one may describe the laws of ${\tilde
P}_{\alpha,\theta,Q}:=P_{\alpha,\theta}\circ \tau_{Q}$ and related
expressions. In view of Remark 1 we will always assume that
$P_{\alpha,\theta}$ is defined by uniform atoms $(U_{i})$, but
allow $\tau_{Q}$ to be based on atoms with a more general
distribution. As mentioned previously there are various techniques
that can be used to identify the laws of random probability
measures. For instance one may calculate its EPPF, identify its
finite dimensional distribution by direct means, or its Laplace
functional. The idea of using Laplace functionals is intuitively
appealing, however it is not particularly suited to handle random
probability measures. It turns out that a more appropriate tool
are generalized Cauchy-Stieltjes transforms~(CS) defined for some
generic random probability measure $\tau$, positive measureable
$g$, positive $z$ and real valued $q$ as
$$
\E[(1+z\tau(g))^{-q}].
$$
Not many results for specific $\tau$ are widely known. Fortunately
there are useful results for the $PD(\alpha,\theta)$ class.
Specifically we shall use the results of Cifarelli and
Regazzini~(1990) and Vershik, Yor and Tsilevich~(2004). Somewhat
less known are the results of James~(2002) who obtains specific
transforms for $PK(\rho)$ and $PK(\rho,\gamma_{\theta})$ models.
These results are extended to larger classes, in a manuscript in
preparation of James, Lijoi and Pr\"unster~(2005). We first describe
what is known for the $PD(\alpha,\theta)$ models. We then describe
the results for the general $PK(\rho,\gamma_{\theta})$ for
$\theta>-\alpha$. This large class of models turn out to be
particularly well suited for coagulation with $PD(\alpha,\theta).$
\subsection{CS for $PD(\alpha,\theta)$}
Note, as can be seen from Remark 1, that given $\tau_{Q}$,
${\tilde P}_{\alpha,\theta,Q}$ is a $PD(\alpha, \theta)$ model
with $H=\tau_{Q}$ fixed. Probably the most widely known CS result
is for the Dirichlet process with shape $\theta H$ where it was
shown by Cifarelli and Regazzini~(1990) that quite remarkably \Eq
\E\[{(1+zP_{0,\theta}(g))}^{-\theta}|H\]={\mbox
e}^{-\int_{\Xcr}d_{\theta}(zg(x))H(dx)} \label{cifreg},\EndEq
where $d_{\theta}$ is defined in~\mref{DPpsi}. Now, key to our
exposition is the following elegant result of Vershik, Yor and
Tsilevich~(2004) for the $PD(\alpha, \theta)$ model with the range
$0<\alpha<1$, $\theta\neq 0$ and otherwise $\theta>-\alpha$ we
have for any $H$ that \Eq
\E[{(1+zP_{\alpha,\theta}(g))}^{-\theta}|H]={\[\int_{\Xcr}{(1+zg(x))}^{\alpha}H(dx)\]}^{-\frac{\theta}{\alpha}}
.\label{keypd} \EndEq To complete the picture for the
$PD(\alpha,\theta)$, a result for the $PD(\alpha,0)$ may be read
from proposition~6.2 of James~(2002) with $n=1$ as, \Eq
\E[{(1+zP_{\alpha,0}(g))}^{-1}|H]=
\frac{\int_{\Xcr}{(1+zg(x))}^{\alpha-1}H(dx)}
{\int_{\Xcr}{(1+zg(x))}^{\alpha}H(dx)}. \label{keystable} \EndEq

\subsection{CS for $PK(\rho,\gamma_{\theta})$ models}
Proposition 6.1 of James~(2002) shows that
$PK(\rho,\gamma_{\theta})$ models for $\theta\neq 0$, $\theta>-1$
have the transform,

\Eq \E[{(1+zP_{K,\theta}(g))}^{-\theta}|H]=
\frac{m_{\theta}(\rho)}{\Gamma(\theta)}\int_{0}^{\infty}{\mbox
e}^{-\int_{\Xcr}\psi(y[1+zg(x)])H(dx)}y^{\theta-1}dy .\label{keyPK}
\EndEq This result generalizes~\mref{keypd}.

\Remark  The result~\mref{keyPK} is actually stated in James~(2002)
for the range $\theta>0$. As the result arises from an identity due
to the gamma function, it extends to the negative range by the same
argument as noted on p. 2309 ``\emph{added in translation}" of
Vershik, Yor and Tsilevich~(2004). Otherwise, one can see this by
first writing,
$$
{(1+zP_{K,\theta}(g))}^{-\theta}={(1+zP_{K,\theta}(g))}^{-(1+\theta)}{(1+zP_{K,\theta}(g))}.
$$
Now noting that $1+\theta>0$ for $\theta>-1$, the gamma identity
applies with $\theta+1$, one then argues as in James~(2002) and
concludes the result by integration by parts.  \EndRemark

\subsection{Generic CS for ${\tilde P}_{\alpha,\theta,Q}$}
Setting $H=\tau_{Q}$ in~\mref{keypd}, we obtain the key formula
for the range $0<\alpha<1$, $\theta\neq 0$ and otherwise
$\theta>-\alpha$
$$
\E[{(1+z{\tilde
P}_{\alpha,\theta,Q}(g))}^{-\theta}|\tau_{Q}]={\[\int_{\Xcr}{(1+zg(x))}^{\alpha}\tau_{Q}(dx)\]}^{-\frac{\theta}{\alpha}}.
$$
It then follows that \Eq \E[{(1+z{\tilde
P}_{\alpha,\theta,Q}(g))}^{-\theta}]=\E[{(1+{\tau_{Q}}(g_{\alpha}))}^{-\frac{\theta}{\alpha}}]
\label{master} \EndEq where
$g_{\alpha}(x)={(1+zg(x))}^{\alpha}-1.$ That is to say, the
peculiar nature of the $PD(\alpha,\theta)$ model for the range
$0<\alpha<1$, $\theta\neq 0$, $\theta>-\alpha$ yields basically a
double Cauchy-Stietltjes formula. Based on the results given in
the previous section an alternative proof for the Theorem~12 in
Pitman~(1999) is almost completely evident. We will describe the
details in the next section. We will then show why the choice of
$PK(\rho,\gamma_{\theta})$ as $Q$ models is quite desirable.
\Remark As we noted earlier explicit CS formula  for more general
$\tau_{Q}$ may be obtained from James~(2002) and James, Lijoi and
Pr\"unster(2005). We note also that the formula for the Dirichlet
model, $PD(0,\theta)$ in~\mref{cifreg} suggests that one needs to
calculate the Laplace functional of a $Q$ which is generally hard.
This makes it more difficult to obtain explicit results for
$PD(0,\theta)$ coagulated by some $Q$. We note from the results in
James~(2002) and James, Lijoi and Pr\"unster~(2005) that the
appearance of Laplace functional calculations is unfortunately
true for many possible candidates as replacements for
$PD(\alpha,\theta).$ In contrast the stable $PD(\alpha,0)$ case is
exceptionally easy. We shall discuss the $PD(0,\theta)$ and
$PD(\alpha,0)$ separately from the other $PD(\alpha,\theta)$
models. \EndRemark

\section{Proof of Pitman's Diagram}
We shall now illustrate how our framework easily yields the
diagram~\mref{Pitman}. First we address the Dirichlet case.
\subsection{$PD(\alpha,\theta)$ coagulated by $PD(0,\theta/\alpha)$}
First apply~\mref{master} with $\tau_{Q}=P_{0,\theta/\alpha}$, then
it is immediate that the formula~\mref{cifreg} applies with
$\theta/\alpha$ in place of $\theta$. To complete the result notice
that
$$
d_{\theta/\alpha}(g_{\alpha}(x))={\theta/\alpha}\ln{(1+zg(x))}^{\alpha}=d_{\theta}(zg(x)).
$$
That is to say the CS transform of order $\theta$ of
$P_{\alpha,\theta}\circ P_{0,\theta / \alpha}$ is given
by~\mref{cifreg} identifying it as a $PD(0,\theta)$ model. The
diagram~\mref{Pitman} is completed by calculating the now obvious
3 EPPF's and applying Bayes rule.
\subsection{$PD(\alpha,\theta)$ coagulated by $PD(\beta,\theta/\alpha)$, $\beta>0$}
For this case apply~\mref{master} with
$\tau_{Q}=P_{\beta,\theta/\alpha}$ for $0<\beta<1$, now
apply~\mref{keypd} to conclude that, in this case, ~\mref{master}
is equivalent to,
$$\E[{(1+{P_{\beta,\theta/\alpha}}(g_{\alpha}))}^{-\frac{\theta}{\alpha}}|H]
={\[\int_{\Xcr}{(1+g_{\alpha}(x))}^{\beta}H(dx)\]}^{-\frac{\theta}{\beta\alpha}}.$$
The result is completed by noting that
$$
{(1+g_{\alpha}(x))}^{\beta}={(1+zg(x))}^{\alpha\beta}.
$$
Hence the CS transform of order $\theta$ of $P_{\alpha,\theta}\circ
P_{\beta,\theta / \alpha}$ is given by the expressions above and now
comparing with \mref{keypd} identifies its law as
$PD(\alpha\beta,\theta)$ model. \Remark The $P_{\alpha,0}\circ
P_{\beta,0}$ case can certainly be obtained easily
using~\mref{keystable} twice. However, we do not believe that there
is a simpler proof than that exhibited in Pitman~(2005). We shall
return to this later. \EndRemark
\section{$PD(\alpha,\theta)$ Coagulated By $PK(\rho,\gamma_{\theta/\alpha})$}
We now obtain a new result as follows. First denote a class of
species sampling random probability measures on $\Xcr$ as \Eq
S_{\alpha,\theta}(\cdot)=\frac{L_{\alpha,\theta}(\cdot)}{T_{L_{\alpha,\theta}}}:=\sum_{i=1}^{\infty}W_{i}\delta_{Z_{i}}(\cdot)
\label{PSrm} \EndEq with
$T_{L_{\alpha,\theta}}=L_{\alpha,\theta}(\Xcr)$ and where
 $(W_{i})\in {\mathcal P}^{\downarrow}_{1}$ and $(Z_{i})$ are iid $H$. Similar to~\mref{absPD} and~\mref{absPK}, the law of
$L_{\alpha,\theta}$, and hence that of $(W_{i})$, is determined by a
power tempered probability measure and satisfies the following
absolute continuity relationship, for every measureable function
$h$, \Eq
\E[h(L_{\alpha,\theta})|H]=c_{\alpha,\theta}m_{\frac{\theta}{\alpha}}(\rho)\E[T^{-\theta}_{L_{\alpha,0}}h(L_{\alpha,0})|H]
\label{absPS} \EndEq where the Laplace functional of the random
measure $L_{\alpha,0}$ is specified by
$$
\E[{\mbox e}^{-L_{\alpha,0}(g)}|H]={\mbox
e}^{-\int_{\Xcr}\tilde{\psi}_{\alpha}(g(x))H(dx)}
$$
with \Eq \tilde{\psi}_{\alpha}(g(x))=\psi({[g(x)]}^{\alpha}).
\label{psia}\EndEq In other words
$L_{\alpha,0}=\mu_{\alpha,0}\circ\mu$, where $\mu_{\alpha,0}$ is a
stable completely random finite measure on $[0,1]$ with index
$0<\alpha<1$ with atoms given by the sequence $(U_{i}),$ and the law
of $\mu$ on $\Xcr$ is specified by~\mref{laplacem}. In particular,
$$
T_{L_{\alpha,0}}\overset {d}=T_{\alpha,0}T^{\frac{1}{\alpha}}
$$
where $T_{\alpha,0}$ is independent of $T$. The density of
$T_{L_{\alpha,0}}$ can be expressed as \Eq f_{T_{L_{\alpha,0}}}(y)=
\int_{0}^{\infty}f_{\alpha}(yt^{-\frac{1}{\alpha}})t^{-\frac{1}{\alpha}}
f_{T}(t)dt=\alpha\int_{0}^{\infty}f_{T}({(y/s)}^{\alpha})s^{-\alpha}y^{\alpha-1}f_{\alpha}(s)ds.
\label{PSden}\EndEq As a by-product, we obtain
$$
c_{\alpha,\theta}m_{\frac{\theta}{\alpha}}(\rho)=1/{\E[T^{-\theta}_{L_{\alpha,0}}]}.
$$
Call the family of laws associated with
 $S_{\alpha,\theta}$, or more specifically the $(W_{i})$, as
 $PS({\rho,\alpha,\theta}).$ In view of~\mref{psia}, the EPPF of the
 $PS({\rho,\alpha,\theta})$ model can be expressed as
\Eq p_{S,\alpha,\theta}(b_{1},\ldots,b_{k})
:=\frac{{(-1)}^{n-k}c_{\alpha,\theta}m_{\frac{\theta}{\alpha}}(\rho)}{\Gamma(\theta+n)}\int_{0}^{\infty}\lambda^{\theta+n-1}{\mbox
e}^{-\psi(\lambda^{\alpha})}\prod_{i=1}^{k}{\tilde
\psi}_{\alpha,b_{i}}(\lambda)d\lambda \label{EPPFPS}\EndEq where
$k$ is the number of blocks formed by the integers
$\{1,\ldots,n\}$ and $\{b_1,\ldots,b_k\}$ are the sizes of the
blocks $\{B_1,\ldots,B_k\}$, respectively. Additionally, similar
to \mref{psim}, for $m=1,\ldots,n,$ $${\tilde
\psi}_{\alpha,m}(\lambda):=\frac{d^m}{d\lambda^m}\psi(\lambda^{\alpha})
:= (-1)^{m-1}\kappa_{\alpha,m}(\lambda), $$ where $
\kappa_{\alpha,m}(\lambda) $ represents the $m$-th cumulant of a
random variable with \emph{tilted} density
$$
{\mbox e}^{{\tilde \psi}_{\alpha}(\lambda)} {\mbox e}^{-\lambda
y}f_{T_{L_{\alpha,0}}}(y),$$ specified by~\mref{PSden}. The
cumulants can be expressed in terms of the moments of this density
$$
\int_{0}^{\infty}y^{m}{\mbox e}^{{\tilde \psi}_{\alpha}(\lambda)}
{\mbox e}^{-\lambda y}f_{T_{L_{\alpha,0}}}(y)dy
$$
using Theile's recursion.\Remark Note that
$PS(\rho_{\beta,0},\alpha,\theta):=PD(\alpha \beta,\theta).$
 \EndRemark
\begin{thm}\label{thm1}Suppose that for $0<\alpha<1$ and $\theta\neq 0$, $\theta>-\alpha,$
$PD(\alpha,\theta)$ is coagulated by
$PK(\rho,\gamma_{\theta/\alpha}).$ Then the following results hold.
\begin{enumerate}
\item[(i)] $PD(\alpha,\theta)$ coagulated by
$PK(\rho,\gamma_{\theta/\alpha})$ is $PS({\rho,\alpha,\theta}),$
with EPPF specified in~\mref{EPPFPS}.
\item[(ii)] Equivalently, suppose that $P_{K,\theta/\alpha}$ is
defined as in~\mref{PKT}, with law determined by $\Pe(\cdot|\rho
H,\gamma_{\theta/\alpha})$, satisfying ~\mref{absPK}. Then the law
of the composition $P_{\alpha,\theta}\circ P_{K,\theta/\alpha}$, is
equivalent to the law of $S_{\alpha,\theta}$, specified
by~\mref{PSrm} and~\mref{absPS}.
\item[(iii)] Suppose there are $K$ blocks $\{A_1,\ldots,A_K\}$ formed by the
integers $\{1,\ldots,n\}$, each with size $a_i$, and $K \geq k$. The
law of the corresponding $Q-FRAG$ kernel is determined by the
(explicit) EPPF
$$
p(a_{1},\ldots,a_{K})=p_{\alpha,\theta}(a_{1},\ldots,a_{K})\times
\frac{\Gamma(\frac{\theta}{\alpha}+1)\Gamma(\theta+n)}{\Gamma(\frac{\theta}{\alpha}+K)\Gamma(\theta+1)}\frac{\int_{0}^{\infty}\lambda^{\theta/\alpha+K-1}{\mbox
e}^{-\psi(\lambda)}\prod_{i=1}^{k}\kappa_{j_{i}}(\lambda)d\lambda}
{\int_{0}^{\infty}\lambda^{\theta+n-1}{\mbox
e}^{-\psi(\lambda^{\alpha})}\prod_{i=1}^{k}\kappa_{\alpha,b_{i}}(\lambda)d\lambda}
$$
where $j_i,i=1,\ldots,k$~(with $\sum_{i=1}^k j_i = K$), is defined
as $\#\{\ell:A_\ell \subseteq B_i\}$.
\end{enumerate}
\end{thm}
\Proof First apply~\mref{master} with
$\tau_{Q}=P_{K,\theta/\alpha}$ and let $C_{1}$ and $C_{2}$ denote
the appropriate constants. Now apply~\mref{keyPK} to get $$
\E[{(1+P_{K,\theta/\alpha}(g_{\alpha}))}^{-\theta/\alpha}|H]=
C_{1}\int_{0}^{\infty}{\mbox
e}^{-\int_{\Xcr}\psi(y[{(1+zg(x))}^{\alpha}])H(dx)}y^{\theta/\alpha-1}dy.
$$ Now apply the transformation $y=w^{\alpha}$ to get
$$ \E[{(1+P_{K,\theta/\alpha}(g_{\alpha}))}^{-\theta/\alpha}|H]
=C_{2}\int_{0}^{\infty}{\mbox
e}^{-\int_{\Xcr}\psi(w^{\alpha}[{(1+zg(x))}^{\alpha}])H(dx)}w^{\theta-1}dw.
$$ Setting
$\tilde{\psi}_{\alpha}(w[(1+zg(x))])=\psi(w^{\alpha}[{(1+zg(x))}^{\alpha})]$,
we see that the CS transform of order $\theta$ of the composition
has the form in~\mref{keyPK} with ${\tilde {\psi}}$ playing the
role of $\psi$. This concludes the result.\EndProof

\subsection{$PD(\alpha,\theta)$ coagulated by $PD(0,\nu)$, Beta Gamma and power tempered normalized Linnik processes}
\label{sec:51} One might be somewhat surprised that
Theorem~\ref{thm1} contains results for $PD(\alpha,\theta)$
coagulated by $PD(0,\nu)$ when $\nu>\theta/\alpha$. In other words,
certain Dirichlet processes are $PK(\rho,\gamma_{\theta/\alpha})$
models. To see this, we recall the Beta Gamma process representation
of Dirichlet processes given in James~(2005). Specifically, if one
chooses a parameter $\nu>\theta/\alpha$, then the law of $\mu$ given
by
$$
\frac{\Gamma(\nu)}{\Gamma(\nu-\theta/\alpha)}T^{-\theta/\alpha}\Pe(d\mu|\rho_{0,\nu}H)
$$
is well defined. Relative to this law, $\mu$ is a Beta Gamma process
with parameters $(\nu H, \theta/\alpha)$ as defined in James~(2005).
Setting $\theta=0$ yields the law of a Gamma process with shape
$\nu$. In any case, James~(2005) shows that normalizing a Beta Gamma
process of this type by its total mass yields a Dirichlet process
with shape $\nu H$ for every $\nu>\theta/\alpha$. Hence the
$PD(0,\nu)$ models are $PK(\rho_{0,\nu},\theta/\alpha)$ models for
$\nu>\theta/\alpha$. Note this equivalence does not hold for
$\nu=\theta/\alpha.$

The corresponding $S_{\alpha,\theta}$ process is obtained by power
tempering of a normalized process, where the law of the process is
determined by
$$
d_{\nu}(\lambda^{\alpha})=\nu\ln(1+\lambda^{\alpha}).
$$
Now we may arrange to have a further scaling which results in the
case where the law of $S_{\alpha,\theta}$ is equivalently obtained
by the power tempering of a normalized Linnik process
subordinator~[see for instance Huillet~(2000, 2003)], where its law
is determined by
$$
\tilde {\psi}_{\alpha}(\lambda):=
d_{\nu}(\lambda^{\alpha}/\nu)=\nu\ln(1+\lambda^{\alpha}/\nu)=\int_{0}^{\infty}(1-{\mbox
 e}^{-\lambda s})l_{\nu,\alpha}(s)ds,
$$
and where
$$
l_{\nu,\alpha}(s)=\frac{\alpha\nu}{s}\phi_{\alpha}(\nu s^{\alpha})
$$
is the L\'evy density of the Linnik process. Specifically, \Eq
\phi_{\alpha}(q)=\E[{\mbox
e}^{-qT^{-\alpha}_{\alpha,0}}]=\sum_{k=0}^{\infty}\frac{1}{\Gamma(1+k\alpha)}{(-q)}^{k}
\label{Mittag}\EndEq is the Mittag-Leffler function or
equivalently the Laplace transform of the random random
$T^{-\alpha}_{\alpha,0}$, where, as before, $T_{\alpha,0}$ is a
stable random variable of index $0<\alpha<1.$ In other words,
$PS(\rho_{0,\nu},\alpha,\theta)=PK(l_{\nu,\alpha},\gamma_{\theta})$
model.

\subsubsection{EPPF calculations}\label{sec:511}
 Note that the above information allows for several descriptions
 of the EPPF of the $PK(l_{\nu,\alpha},\gamma_{\theta})$ model and
 hence the corresponding $Q-FRAG$. First using~\mref{Mittag} one
 has that
the $m$-th cumulant of an exponentially tilted Linnik random
variable can be expressed as
$$
\kappa_{\alpha,m}(u)={\alpha\nu}\int_{0}^{\infty}s^{m-1}{\mbox
e}^{-us}\phi_{\alpha}(\nu
s^{\alpha})ds={\alpha\nu}u^{-m}\sum_{l=0}^{\infty}\frac{\Gamma(m+l\alpha)}{\Gamma(1+l\alpha)}u^{-l\alpha}
{(-\nu)}^{l}
$$
or
$$
\kappa_{\alpha,m}(u)=\alpha\nu
u^{-m}\int_{0}^{\infty}\int_{0}^{\infty}s^{m-1}{\mbox e}^{-s-\nu
{(s/u)}^{\alpha}t^{-\alpha}}f_{\alpha}(t)dtds.
$$
The general EPPF of the $PK(l_{\nu,\alpha},\gamma_{\theta})$ model
can be written as
\begin{equation}\label{LEPPF}
\frac{c_{\alpha,\theta}\Gamma(\nu)}{\Gamma(\nu-\theta/\alpha)}\frac{\nu^{\theta/\alpha+k}\alpha^{k-1}}{\Gamma(\theta+n)}\int_{0}^{\infty}y^{\theta/\alpha-1}{(1+y)}^{-\nu}
\prod_{j=1}^{k}R_{n_{j}}(y|\alpha)dy,
\end{equation}
where
$$
R_{n_{j}}(y|\alpha)=\sum_{l=0}^{\infty}\frac{\Gamma(n_{j}+l\alpha)}{\Gamma(1+l\alpha)}
{(-y)}^{-l}= \int_{0}^{\infty}\int_{0}^{\infty}s^{n_{j}-1}{\mbox
e}^{-s- y^{-1}{s}^{\alpha}t^{-\alpha}}f_{\alpha}(t)dtds.
$$
Furthermore in the case of $\alpha=1/2$, we may follow the results
for Brownian excursion in Pitman~(2003, Section 8), to obtain
$$
R_{n_{j}}(y|1/2)=\frac{\Gamma(2n_{j})2^{-n_{j}+1/2}}{\Gamma(1/2)}\int_{0}^{\infty}h_{-2n_{j}}(x/y){\mbox
e}^{-x^{2}/2}dx
$$
where
$$
h_{-2n_{j}}(x)=\frac{2^{n_{j}-1}}{\Gamma(2n_{j})}\int_{0}^{\infty}s^{n_{j}-1}{\mbox
e}^{-s- x\sqrt{2s}}ds
$$
is a Hermite function.

\begin{prop}\label{propo}Suppose that for $0<\alpha<1,$ $PD(\alpha,\theta)$
is coagulated by $PD(0,\nu)$ for $\nu>\theta/\alpha$, then the
following results hold
\begin{enumerate}
\item[(i)] $PD(\alpha,\theta)$ coagulated by
$PD(0,\nu)$ is $PK(l_{\nu,\alpha},\gamma_{\theta})$,  with EPPF
specified in~\mref{LEPPF}.
\item[(ii)] Suppose there are $K$ blocks $\{A_1,\ldots,A_K\}$ formed by the
integers $\{1,\ldots,n\}$, each with size $a_i$, and $K \geq k$.
The law of the corresponding $Q-FRAG$ kernel is determined by the
EPPF \Eq p_{\alpha,\theta}(a_{1},\ldots,a_{K})\times
\frac{\Gamma(\nu-\theta/\alpha)\Gamma(\theta+n)}{c_{\alpha,\theta}\nu^{\theta/\alpha}\alpha^{k-1}\Gamma(\nu+K)}\frac{\prod_{i=1}^k(j_i-1)!}
{\int_{0}^{\infty}y^{\theta/\alpha-1}{(1+y)}^{-\nu}
\prod_{j=1}^{k}R_{b_{j}}(y|\alpha)dy} \label{DPEPPF}\EndEq where
$j_i,i=1,\ldots,k$~(with $\sum_{i=1}^k j_i = K$), is defined as
$\#\{\ell:A_\ell \subseteq B_i\}$.
\end{enumerate}
\end{prop}

\Remark It is evident from Proposition~\ref{propo} that it is not
easy to obtain the denominator in~\mref{DPEPPF} by summing out
appropriately over the numerator. This is despite the fact that
both the EPPF's in the numerator have nice Gibbs form. Hence again
this points to the difficulties of a direct combinatorial
argument. On the other hand, the results establish some rather
peculiar combinatorial identities. \EndRemark

\section{$PD(\alpha,0)$ coagulated by
$PK(\rho)$} In view of the arguments in Bertoin and LeGall~(2003)
and Pitman~(2005, p. 115) concerning the Bolthausen-Sznitman~(1998)
coalescent, that is a description of $PD(\alpha,0)$ coagulated by
$PD(\beta,0)$, it is easy to extend this to the case of
$PD(\alpha,0)$ coagulated by $PK(\rho)$

\begin{thm}Suppose that for $0<\alpha<1$
$PD(\alpha,0)$ is coagulated by $PK(\rho).$ Then the diagram,
according to~\mref{arrow3}, of this process is described by setting
$\theta=0$ in Theorem~\ref{thm1}.
\end{thm}

\Proof The proof proceeds along the same lines as Pitman~(2005, p.
115) and Bertoin and Le Gall~(2003, p. 272). That is,
$$
\frac{\mu_{\alpha,0}(P_{K}(\cdot))}{T_{\alpha,0}}\overset
{d}=\frac{\mu_{\alpha,0}(\mu(\cdot))}{\mu_{\alpha,0}(T)}\overset{d}=
\frac{L_{\alpha,0}(\cdot)}{T_{L_{\alpha,0}}}:=S_{\alpha,0}(\cdot)
$$
where $L_{\alpha,0}$, $S_{\alpha,0}$ are as described in the
beginning of this section. \EndProof

\Remark The description of $PD(\alpha,0)$ coagulated by
$PD(0,\nu)$ for $\nu>0$ is obtained from Proposition~\ref{propo}
with $\theta=0.$\EndRemark

\section{Some comments about $PD(0,\theta)$ coagulated by general $Q$}
It was noted earlier that obtaining results for the $PD(0,\theta)$
coagulated by some $Q$ does not readily follow from a CS type
analysis. We will now briefly describe how one can obtain the finite
dimensional distributions of such compositions. Note again that a
Dirichlet process coagulated by a $Q$, i.e., $P_{0,\theta}\circ
\tau_{Q}$ given $\tau_{Q}$, is a Dirichlet process with shape
$\theta \tau_{Q}$  on $\Xcr.$ Hence it follows that for any
measureable partition $C_{1},\ldots, C_{m}$ of $\Xcr$ the finite
dimensional distribution of $P_{0,\theta}\circ \tau_{Q}$ given
$\tau_{Q}$ is specified by the joint Dirichlet density of
$Y_{i}=P_{0,\theta}(\tau_{Q}(C_{i}))$ for $i=1,\ldots,m$, which is
given by
$$
f(y_{1},\ldots,y_{m}|\tau_{Q})=\frac{\Gamma(\theta)}{\prod_{i=1}^{m}\Gamma(\theta
z_{i})}\prod_{i=1}^{m}y^{\theta z_{i}-1}_{i}
$$
where $z_{i}=\tau_{Q}(C_{i}),$ and $(Y_{1},\ldots,Y_{m})\in
\Scr_{m}=\{(a_{i})_{i\leq m}:0<a_{i}<1,\sum_{i=1}^{m}a_{i}=1\}.$
This leads to a general description of the finite dimensional
distributions.

\begin{prop}\label{prop2}Suppose that $P_{0,\theta}$ denotes a Dirichlet Process on
$[0, 1]$ with shape $\theta U$ and $U$ is a uniform distribution.
Suppose further that $\tau_{Q}$ is a random probability measure
on~$\Xcr$. Then, for a measureable partition $C_{1},\ldots, C_{m}$
of $\Xcr$, the distribution of $P_{\theta, Q}=P_{0,\theta}\circ
\tau_{Q}$ is specified by its finite-dimensional distribution $$
f_{\theta,Q}(y_{1},\ldots,y_{m})=\int_{\Scr_{m}
}f_{Q}(z_{1},\ldots,z_{m})\frac{\Gamma(\theta)}{\prod_{i=1}^{m}\Gamma(\theta
z_{i})}\prod_{i=1}^{m}y^{\theta z_{i}-1}_{i}dz_{i}$$ where
$Y_{i}=P_{\theta, Q}(C_{i})$, and $f_{Q}$ denotes the joint density
of $Z_{i}=Q(C_{i})$ for $i=1,\ldots,m$. If $Q$ is a species sampling
model, then this equates to a description of the law of
$PD(0,\theta)$ coagulated by $Q$.\end{prop}

\noindent Naturally the utility of this result requires knowledge of
the finite-dimensional distribution of $Q$. Below we describe two
special cases where $Q=PD(1/2,\eta)$ for $\eta>-1/2$ and
$Q=PD(0,\nu)$ for $\nu>0.$

\begin{prop}
The distribution of $PD(0,\theta)$ coagulated by $PD(1/2,\eta)$
for $\eta>-1/2$  is determined by the distribution of
$P_{0,\theta}\circ P_{1/2,\eta}$ with finite dimensional
distribution,
$$\frac{\[\prod_{1=1}^{m}p_i\]\Gamma(\eta+m/2)}{\pi^{(m-1)/2}\Gamma(\eta+1/2)}
\int_{\Scr_{m} }\frac{\Gamma(\theta)}{\prod_{i=1}^{m}\Gamma(\theta
z_{i})}\prod_{i=1}^{m}y^{\theta z_{i}-1}_{i}
\frac{\prod_{i=1}^{m}z^{-3/2}_{i}} {{(p^{2}_{1}/z_{1}+\cdots
+p^{2}_{m}/z_m)}^{\eta+m/2}} dz_{i},$$ where $p_{i}=H(C_{i})$ and
otherwise the notation is as in Proposition~\ref{prop2}.\end{prop}

\Proof The distribution of $(Q(C_{i}))$ when $Q$ is $PD(1/2,\eta)$
is given by Theorem~3.1 of Carlton~(2002). \EndProof

We now describe the Dirichlet case.

\begin{prop}
The distribution of $PD(0,\theta)$ coagulated by $PD(0,\nu)$   is
determined by the distribution of $P_{0,\theta}\circ P_{0,\nu}$
with finite dimensional distribution, $$ \int_{\Scr_{m}
}\frac{\Gamma(\theta)}{\prod_{i=1}^{m}\Gamma(\theta
z_{i})}\prod_{i=1}^{m}y^{\theta
z_{i}-1}_{i}\frac{\Gamma(\nu)}{\prod_{i=1}^{m}\Gamma(\nu
p_{i})}\prod_{i=1}^{m}z^{\nu p_{i}-1}_{i} dz_{i}$$ where
$p_{i}=H(C_{i})$ and otherwise the notation is as in
Proposition~\ref{prop2}.\end{prop}

\Remark It is interesting to note that the dynamics of the
coagulation of two or more Dirichlet processes may also be
explained, in perhaps a more informative way, as a Chinese
restaurant franchise process of Teh, Jordan, Beal and Blei~(2006)
when there is one franchise. In their setup, one has $F|\tau_{Q}$ is
Dirichlet process with shape $(\theta \tau_{Q})$ and $\tau_{Q}$ is a
Dirichlet process with shape, say, $\eta H.$ The distribution of $F$
is characterized via a Chinese restaurant franchise with one
franchise. It follows from our observations that $F\overset
{d}=P_{0,\theta}\circ P_{0,\nu},$ leading to an equivalence. Based
on this observation, all the processes discussed here lead to some
type of Chinese restaurant franchise process, and therefore have
potential applications in machine learning and related areas.
\EndRemark

\vskip0.2in \centerline{\Heading References} \vskip0.2in \tenrm
\def\smc{\tensmc}
\def\sl{\tensl}
\def\bf{\tenbold}
\baselineskip0.15in

\Ref \by    Aldous, D. J. and Pitman, J.  \yr    1998 \paper The
standard additive coalescent \jour \AnnProb \vol   26 \pages
1703-1726 \EndRef

\Ref \by    Bertoin, J.  \yr    2002 \paper Self-similar
fragmentations \jour Ann. Inst. H. Poincaré Probab. Statist. \vol 3
\pages 319-340 \EndRef

\Ref \by    Bertoin, J. and Le Gall, J.-F.  \yr    2003 \paper
Stochastic flows associated to coalescent processes \jour Probab.
Theory Related Fields \vol   126 \pages 261-288\EndRef

\Ref \by Bertoin, J. and Le Gall, J.-F.  \yr    2005 \paper
Stochastic flows associated to coalescent processes. II. Stochastic
differential equations \jour Ann. Inst. H. Poincaré Probab. Statist.
\vol   3 \pages 307-333\EndRef

\Ref \by Bertoin, J. and Goldschmidt, C.  \yr    2004 \paper Dual
random fragmentation and coagulation and an application to the
genealogy of Yule processes. In \textit{Mathematics and computer
science III: Algorithms, Trees, Combinatorics and Probabilities}, M.
Drmota, P. Flajolet, D. Gardy, B. Gittenberger (editors), pp.
295-308. Trends Math., Birkh\"auser, Basel  \EndRef

\Ref \by Bertoin, J. and Pitman, J.  \yr    2000 \paper Two
coalescents derived from the ranges of stable subordinators \jour
Electron. J. Probab. \vol   7 \pages 1-17\EndRef

\Ref \by Bolthausen, E. and Sznitman, A.-S.  \yr    1998 \paper On
Ruelle's probability cascades and an abstract cavity method \jour
Comm. Math. Phys. \vol   197 \pages 247-276\EndRef

\Ref \by Carlton, M. A.  \yr    2002 \paper A family of densities
derived from the three-parameter Dirichlet process \jour J. Appl.
Probab. \vol 39 \pages 764-774\EndRef

\Ref \by    Cifarelli, D. M. and Regazzini, E. \yr    1990 \paper
Distribution functions of means of a Dirichlet process \jour
\AnnStat \vol   18 \pages 429-442 \EndRef

\Ref \by Dong, R., Martin, J. and Goldschmidt, C.  \yr    2005
\paper Coagulation-fragmentation duality, Poisson-Dirichlet
distributions and random recursive trees, arXiv math.PR/0507591,
2005. \\ Available at http://arxiv.org/abs/math.PR/0507591. \EndRef

\Ref \by Huillet, T.  \yr 2000 \paper On Linnik's continuous-time
random walks \jour J. Phys. A \vol 33 \pages 2631-2652\EndRef

\Ref \by Huillet, T.  \yr 2003 \paper Energy cascades as branching
processes with emphasis on Neveu's approach to Derrida's random
energy model \jour Adv. in Appl. Probab. \vol 35 \pages
477-503\EndRef

\Ref \by     James, L.F. \yr     2002 \paper Poisson process
partition calculus with applications to exchangeable models and
Bayesian nonparametrics,
arXiv:math.PR/0205093, 2002.\\
Available at \texttt{http://arxiv.org/abs/math.PR/0205093} \EndRef

\Ref \by James, L.F.  \yr 2005 \paper Functionals of Dirichlet
processes, the Cifarelli-Regazzini identity and Beta-Gamma processes
\jour \AnnStat \vol 33 \pages 647-660\EndRef

\Ref \by     James, L.F., Lijoi, A. and Pr\"unster, I. \yr 2005
\paper Cifarelli-Regazzini/Markov-Krein type identities and
distributional results for functionals of normalized random
measures. Manuscript in preparation \EndRef

\Ref \by Kingman, J. F. C. \yr 1975 \paper Random discrete
distributions \jour J. R. Stat. Soc. Ser. B \vol 37 \pages
1-22\EndRef

\Ref \by    Pitman, J. \yr 1996 \paper Some developments of the
Blackwell-MacQueen urn scheme. In \textit{Statistics, Probability
and Game Theory}, T.S. Ferguson, L.S. Shapley and J.B. MacQueen
(editors), IMS Lecture Notes-Monograph series, Vol. 30, pp. 245-267,
Inst. Math. Statist., Hayward, CA \EndRef

\Ref \by    Pitman, J.  \yr    1999 \paper Coalescents with multiple
collision \jour \AnnProb \vol   27 \pages 1870-1902 \EndRef

\Ref \by Pitman, J. \yr    2003 \paper Poisson-Kingman partitions.
In \textit{Statistics and science: a Festschrift for Terry Speed},
D.R. Goldstein (editor), IMS Lecture Notes-Monograph series, Vol.
40, pp. 1-34, Inst. Math. Statist., Hayward, CA \EndRef

\Ref \by    Pitman, J. \yr    2005 \paper Combinatorial Stochastic
Processes. Lecture Notes in Mathematics and in Probability Surveys,
Springer. \\
Available at
\texttt{http://bibserver.berkeley.edu/csp/april05/bookcsp.pdf}
\EndRef

\Ref \by    Pitman, J. and Yor, M. \yr    1997 \paper The
two-parameter Poisson-Dirichlet distribution derived from a stable
subordinator \jour \AnnProb \vol   25 \pages 855-900 \EndRef

\Ref \by    Schweinsberg, J. \yr 2000 \paper Coalescents with
simultaneous multiple collisions \jour Electron. J. Probab. \vol 5
\pages 1-50 \EndRef

\Ref \by    Teh, Y.W., Jordan, M.I., Beal, M.J. and Blei, D.M. \yr
2006 \paper Hierarchical Dirichlet Processes. Available at
\texttt{http://stat-www.berkeley.edu/tech-reports/index.html}\\To
appear in \textit{J. Amer. Statist. Assoc} \EndRef

\Ref \by    Vershik, A.M., Yor, M. and Tsilevich, N.V.  \yr    2004
\paper On the Markov-Krein identity and quasi-invariance of the
gamma process \jour J. Math. Sci. \vol   121 \pages 2303-2310
\EndRef

\vskip0.75in

\smc

\Tabular{ll}

Man-Wai Ho\\
Department of Statistics and Applied Probability \\
National University of Singapore\\
6 Science Drive 2\\
Singapore 117546\\
Republic of Singapore\\
\rm stahmw\at nus.edu.sg\\

\\
\EndTabular

\Tabular{ll}

Lancelot F. James\\
The Hong Kong University of Science and Technology\\
Department of Information and Systems Management\\
Clear Water Bay, Kowloon\\
Hong Kong\\
\rm lancelot\at ust.hk\\

\\
\EndTabular

\Tabular{ll}

John W. Lau\\
Department of Mathematics\\
University of Bristol Bristol\\
BS8 1TW\\
United Kingdom\\
\rm John.Lau\at bristol.ac.uk\\

\EndTabular

\end{document}